\documentstyle{article}\begin{document}\centerline{\bf  Generalized  Bessel recursion relations }\vskip .3in

\centerline{M.L. Glasser}\vskip .4in

\centerline{ Department of Physics}

\centerline{Clarkson University}

\centerline{Potsdam, NY (USA)}\vskip .5in

\centerline{\bf ABSTRACT}

\vskip .1in
\begin{quote}

This paper presents the equality of finite index sums of Bessel func- tions containing arbitrary numbers of terms. These reduce to the familiar three term recursion formulas in simple cases.
\end{quote}



\newpage

\section{Introduction}
The  motivation for this note was the observation that the basic recursion relation for the modified Bessel function $K$[1],  
$$K_0(z)+\left(\frac{2}{z}\right)K_1(z)=K_2(z)$$
can be expressed as the symmetry with respect to $m=0$ and $n=1$ of the sum
$$\sum_{k=0}^n K_{k-m-1}(z)\left(\frac{z}{2}\right)^{k+m}.\eqno(1)$$

The attempt to generalize this to arbitrary $m$ and $n$ led to our principal result\vskip .1in

\noindent
{\bf Theorem 1}\vskip .1in
For positive integers $m$ and $n$ the expression
$$(n+1)!\sum_{k=0}^n\frac{1}{k!}{m+k+1\choose{m}}K_{k-m-1}(z)\left(\frac{z}{2}\right)^{k+m}\eqno(2)$$
is symmetric with respect to $m$ and $n$.\vskip .1in

\noindent
This will be proven in the following section and some similar results presented in the concluding paragraph.

\section{Calculation}

Consider the sum
$$F(n,p,q)=\frac{(n+q+1)!}{q!(q+1)!}\sum_{k=0}^p \frac{(q+k+1)!}{(k+1)!}\frac{(n+k)!}{k!}\eqno(3)$$
for $p,q,n\in {\cal{Z}}^+$. One finds that, e.g.
$$F(1,p,q)=\frac{(p+q+2)!}{p!q!}$$
$$F(2,p,q)=\frac{(p+q+2)!}{p!q!}[6+2(p+q)+pq]$$
and by induction on $n$ one obtains
\vskip .1in

\newpage

\noindent
{\bf Lemma 1}\vskip .1in
$$\frac{p!q!}{(p+q+2)!}F(n,p,q)$$
is a polynomial $P(p,q)=P(q,p)$ of degree $n-1$ in $p$ and $q$.

Next, by interchanging the order of summation and invoking lemma 1, one has
\vskip .1in

\noindent
{\bf Lemma 2}\vskip .1in

$$G(p,q,z)=\sum_{n=0}^{\infty}\frac{1}{(n!)^2}F(n,p,q)z^n=\sum_{k=0}^p {q+k+1\choose{q}}\;_2F_1(k+1,q+2;1;z)$$
is analytic for $ |z|<1$ and symmetric with respect to $p$ and $q$.

\vskip .1in
Finally, noting that[2]
$$\int_0^{\infty}J_0(z\sqrt{x})\;_2F_1(k+1,q+2;1;-x)dx=\frac{2^{-k-q} z^{k+q+1}}{k!(q+1)!}K_{k-q-1}(z)\eqno(4)$$
(changing $q$ to $m$ and $p$ to $n$) we have Theorem 1.
 
 For example, with $m=0$ we get the possibly new summation
 $$\sum_{k=0}^n\frac{1}{k!}K_{k-1}(z)(z/2)^k=\frac{1}{n!}K_{n+1}(z)(z/2)^n.\eqno(5)$$
 
 Setting $z=-ix$ in the relation
 $$K_{\nu}(z)= \frac{\pi}{2}i^{\nu+1}[J_{\nu}(iz)+iY_{\nu}(iz)]\eqno(6)$$
 after a small manipulation one obtains 
 \vskip .1in
 
 \noindent
 {\bf Theorem 2}\vskip .1in
 
 $$(-1)^m(n+1)!\sum_{k=0}^n\frac{1}{k!}{m+k+1\choose{m}}\, J_{k-m-1}(x)(x/2)^{k+m}\eqno(7)$$
 $$(-1)^m(n+1)!\sum_{k=0}^n\frac{1}{k!} {m+k+1\choose{m}}\, Y_{k-m-1}(x)(x/2)^{k+m}\eqno(8)$$
 are both symmetric  with respect to $m$ and $n$.\vskip .1in
 \newpage
 \noindent
 {\bf Corollary }\vskip .1in
 
 $$\sum_{k=0}^{n}\frac{1}{k!}\, {\cal{C}}_{k-1}(x)(x/2)^k
 =-\frac{1}{n!}\, {\cal{C}}_{n+1}(x)(x/2)^n\eqno(9)$$
 where ${\cal{C}}=aJ +b Y$.
  \vskip .2in
 \section{Discussion}
 
 Analogous sum relations can be obtained by other means. For example,
 let us start with the hypergeometric summation formula[3]
$$\;_3F_2(-n,1,a;3-a,n+3;-1)=\frac{(n+2)n!}{2(a-1)\Gamma(a-2)}\left[\frac{\Gamma(a-1)}{(n+1)!}+(-1)^n\Gamma(a-n-2)\right].\eqno(9)$$
But, 
$$\;_3F_2(-n,1,a;3-a,n+3;-1)=\frac{n!(n+2)!}{\Gamma(a-2)\Gamma(a)}\sum_{k=1}^{n+1} (-1)^{k+1}\frac{\Gamma(a-1+k)\Gamma(a-1-k)}{\Gamma(n+k)\Gamma(n-k)}.\eqno(10)$$
With $n$ replaced by $n-1$ and $a=(s+n)/2+1$, the first term of (9) is half of what would be the $k=0$ term of the sum in (10) and one has

$$\sum_{k=0}^n(-1)^k(2-\delta_{k,0})\frac{\Gamma\left(\frac{s+n}{2}-k\right)\Gamma\left(\frac{s+n}{2}+k\right)}{(n-k)!(n+k)!}=\frac{(-1)^n}{n!}\Gamma\left(\frac{s+n}{2}\right)\Gamma\left(\frac{s-n}{2}\right).\eqno(11)$$
Next we take the inverse Mellin transform of both sides, noting that
$$\int_{c-i\infty}^{c+i\infty}\frac{ds}{2\pi i}(2/x)^s\Gamma\left(\frac{s+n}{2}-k\right)\Gamma\left(\frac{s+n}{2}+k\right)=4x^nK_{2k}(x)\eqno(12)$$
$$\int_{c-i\infty}^{c+i\infty}\frac{ds}{2\pi i}(2/x)^s\Gamma\left(\frac{s-n}{2}\right)\Gamma\left(\frac{s+n}{2}\right)=4K_n(x).\eqno(13)$$
Consequently,
$$K_n(x)=\left(\frac{x}{2}\right)^n\sum_{k=0}^{n}(-1)^{k+n}n!\frac{(2-\delta_{k,0})}{(n-k)!(n+k)!}K_{2k}(x).\eqno(14)$$

 Since many integrals of the Gauss hypergeometric function are known, one of the most extensive tabulations being[2], Lemma 2 is the gateway to a myriad of unexpected finite sum identities involving various classes of special functions. We conclude by listing a small selection..
 
 From[2] 
 
 $$\int_0^{\infty}(1-e^{-t})^{\lambda-1}e^{-xt}\;_2F_1(k+1,m+2;1;ze^{-t})dt$$
 $$=B(x,\lambda)\;_3F_2(k+1,m+2,x;1,x+\lambda; z) \eqno(15)$$
 and one has the symmetry of
$$ \sum_{k=0}^n {m+k+1\choose{m}}\;_3F_2(k+1,m+2,x;1,x+\lambda;z)\eqno(16)$$
For example for $m=0$
$$\sum_{k=0}^n\;_3F_2(k+1,2,x;1,x+\lambda;z)=(n+1)\;_2F_1(n+2,x;x+\lambda;z).\eqno(17)$$

Similarly,
$$\frac{n!(n+1)!}{\Gamma(n+2-a)}\sum_{k=0}^n\frac{(m+k+1)!\Gamma(k+1-a)}{k!(k+1)!}.$$
$$=\frac{m!(m+1)!}{\Gamma(m+2-a)}\sum_{k=0}^m\frac{(n+k+1)!\Gamma(k+1-a)}{k!(k+1)!}\eqno(18)$$

$$\sum_{k=0}^n{m+k+1\choose{m}}=\sum_{k=0}^m{n+k+1\choose{n}}.\eqno(19)$$

$$\sum_{k=0}^n{m+k+1\choose{m}}\;_3F_2(k+1,m+2,a;1,a+b;z)$$
$$=\sum_{k=0}^m{n+k+1\choose{n}}\;_3F_2(k+1,n+2,a;1,a+b;z).\eqno(20)$$

$$\sum_{k=0}^n\;_3F_2(k+1,2,a;1,a+b;1)=\frac{(n+1)\Gamma(b-n-2)\Gamma(a+b)}{\Gamma(a+b-n-2)\Gamma(b)}.\eqno(21)$$

$$\sum_{k=0}^n\frac{(p+k)!}{k!}=\frac{(n+p+1)!}{(p+1)n!},\quad p=0,1,2,\cdots\eqno(22)$$

$$\sum_{k=0}^n{m+k+1\choose{m}}z^{(k+m)/2}S_{-k-m-2,k-m-1}(z)$$
$$=\sum_{k=0}^m{n+k+1\choose{n}}z^{(k+n)/2}S_{-k-n-2,k-n-1}(z).\eqno(23)$$

$$\sum_{k=0}^n{m+k+1\choose{m}}z^{(k+m)/2}W_{-k-m-2,k-m-1}(z)$$
$$=\sum_{k=0}^m{n+k+1\choose{n}}z^{(k+n)/2}W_{-k-n-2,k-n-1}(z).\eqno(24)$$

\section{References}

\noindent
[1] G.E. Andrews, R. Askey and R. Roy, {\it Special Functions} [Cambridge University Press, 1999]

\noindent
[2]  A.P. Prudnikov, Yu. A.Brychkov and O.I. Marichev,{\it Integrals and Series, Vol. 3} [Gordon and Breach, NY 1986] Section 2.21.1.
 
 \noindent
 [3] Ibid. Section (2.4.1).

\end{document}